\documentclass[12pt,a4paper,reqno]{amsart}
\usepackage{amssymb}
\usepackage{dcpic,pictex}
\usepackage{multirow}
\usepackage{graphicx}
\usepackage{slashbox}
\usepackage{cite}

\usepackage{hyperref}
\usepackage[main=english,french]{babel}
\newtheorem{definition}{Definition}[section]
\newtheorem{example}{Example}[section]
\newtheorem{theorem}{Theorem}[section]
\newtheorem{remark}{Remark}[section]
\newcommand{\D}{\mathbb{D}}

\newcommand{\I}{\mathbb{I}}
\newcommand{\R}{\mathbb{R}}
\newcommand{\C}{\mathbb{C}}
\newcommand{\Com}{\mathbb{C}}
\newcommand{\N}{\mathbb{N}}

\begin{document}

%\selectlanguage{english}

\title[General fractional integrals and derivatives]{General fractional integrals and 
        derivatives of arbitrary order}

\author{Yuri Luchko}
%\address{Beuth University of Applied Sciences Berlin}
\curraddr{Beuth Technical University of Applied Sciences Berlin,  
     Department of  Mathematics, Physics, and Chemistry,  
     Luxemburger Str. 10,  
     13353 Berlin,
Germany}
\email{luchko@beuth-hochschule.de}

%    General info
%\UDC{519.21}
%\date{04.07.2020}
\dedicatory{}
% \thanks{Paper accepted for publication in {\it Fract. Calc. Appl. Anal.}, to appear at https://www.degruyter.com/view/j/fca}
% published in 
% Fract. Calc. Appl. Anal.,Vol. 23, No 4 (2020), pp. 939–966, DOI: 10.1515/fca-2020-0049

\begin{abstract}
In this paper, we introduce the general fractional integrals and derivatives of arbitrary order and study some of their basic properties and particular cases. First, a suitable generalization of the Sonine condition is presented and some important classes of the kernels that satisfy this condition are introduced. Whereas the kernels of the general fractional derivatives with these kernels possess the integrable singularities at the point zero, the kernels of the general fractional integrals can be - depending on their order - both singular and continuous at the origin.   For the general fractional integrals and derivatives of arbitrary order with the kernels introduced in this paper, two fundamental theorems of  fractional  calculus are formulated and proved. 
\end{abstract}

\maketitle

%%%%%%%%%%%%%%%%%%%%%%%%%%%%%%%%%%%%%%%%%%%%%%%%
%\begin{document}

\section{Introduction}\label{sec:1}

In his papers \cite{Abel1,Abel2}, Abel derived and studied a mathematical model for the  tautochrone problem in form of the following integral equation (in slightly different notations):
\begin{equation}
\label{A1}
f(t) = \frac{1}{\sqrt{\pi}}\, \int_0^t \frac{\phi^\prime(\tau)\, d\tau}{\sqrt{t-\tau}}.
\end{equation}
In fact, he considered even more general integral equation  
\begin{equation}
\label{A2}
f(t) =\frac{1}{\Gamma(1-\alpha)}\, \int_0^t \frac{\phi^\prime(\tau)\, d\tau}{(t-\tau)^\alpha}
\end{equation}
under an implicit restriction $0<\alpha <1$. It is easy to see that the right-hand side of \eqref{A2} is the operator that is nowadays referred to as the Caputo fractional derivative $ _*D_{0+}^{\alpha}$ of the order $\alpha,\ 0<\alpha <1$. 
Abel's solution formula to the equation \eqref{A2} is nothing else than the operator nowadays  called the Riemann-Liouville fractional integral $I_{0+}^\alpha$ of the order $\alpha >0$: 
\begin{equation}
\label{A4}
\phi(t) = \frac{1}{\Gamma(\alpha)}\, \int_0^t (t-\tau)^{\alpha-1}\, f(\tau)\, d\tau \ =: \ (I_{0+}^\alpha\, f)(t),\ t>0.
\end{equation}
In the modern notations, the formulas \eqref{A2} and \eqref{A4} build the 2nd fundamental theorem of FC for the Caputo fractional derivative of a function taking the value zero at the point zero:
\begin{equation}
\label{A5}
(I_{0+}^\alpha\, f)(t) = (I_{0+}^\alpha\,  _*D_{0+}^{\alpha}\, \phi)(t) = \phi(t) -\phi(0) = \phi(t),
\end{equation}
where validity of the condition $\phi(0) = 0$ follows from construction of Abel's mathematical model for the  tautochrone problem. For more details regarding Abel's results and derivations  presented in \cite{Abel1,Abel2} see the recent paper \cite{Pod}.

To solve  the integral equation \eqref{A2}, in the paper \cite{Abel2} published in 1826, Abel  employed the relation
\begin{equation}
\label{A3}
(h_{\alpha}*h_{1-\alpha})(t) = \{ 1 \},\ t>0,\ h_{\alpha}(t):= \frac{t^{\alpha -1}}{\Gamma(\alpha)},\ \alpha>0,
\end{equation}
where the operation $*$ stands for the Laplace convolution
\begin{equation}
\label{con}
(f\, *\, g)(t) = \int_0^t f(t-\tau)g(\tau)\, d\tau
\end{equation}
and $\{ 1 \}$ is the function that is identically equal to $1$ for $t\ge 0$.

In the paper \cite{Son} published in 1884, Sonine recognized that the relation \eqref{A3} is the most crucial ingredient of  Abel's solution method that can be generalized and applied for analytical treatment of a larger class of integral equations. In place of \eqref{A3}, Sonine  considered a pair of functions $\kappa,\ k$ (Sonine kernels)  that satisfy the relation 
\begin{equation}
\label{Son}
(\kappa \, *\, k )(t) = \{1 \},\ t>0.
\end{equation}
In what follows, we denote the set of the Sonine kernels by $\mathcal{S}$. For a given Sonine kernel $\kappa$, the kernel  $k$ that satisfies the Sonine condition \eqref{Son} is called its associate Sonine kernel. 
Following Abel's solution method, Sonine showed that the integral equation 
\begin{equation}
\label{Son1}
f(t) = \int_0^t \kappa(t-\tau)\phi(\tau)\, d\tau \, = \, (\kappa\, *\, \phi)(t)
\end{equation}
has a solution in the form
\begin{equation}
\label{Son2}
\phi(t) = \frac{d}{dt}\, \int_0^t k(t-\tau)f(\tau)\, d\tau \, = \, \frac{d}{dt} (k\, *\, f)(t),
\end{equation}
provided the kernels $\kappa,\ k$ satisfy the Sonine condition \eqref{Son}. Indeed, we get
$$
(k\, *\, f)(t) = (k\, *\, \kappa\, *\, \phi)(t) = (\{1 \} \, *\, \phi)(t) = \int_0^t \phi(\tau)\, d\tau
$$
that immediately leads to the formula \eqref{Son2}. Of course, any concrete realization of the Sonine schema requires a precise characterization of the Sonine kernels and the spaces of functions where the operators from the right-hand sides of \eqref{Son1} and \eqref{Son2} are well defined. In \cite{Son}, Sonine introduced a large class of the Sonine kernels in the form
\begin{equation}
\label{3-3}
\kappa(t) = h_{\alpha}(t) \cdot \, \kappa_1(t),\ \kappa_1(t)=\sum_{k=0}^{+\infty}\, a_k t^k, \ a_0 \not = 0,\ 0<\alpha <1,
\end{equation}
\begin{equation}
\label{3-4}
k(t) = h_{1-\alpha}(t) \cdot k_1(t),\ k_1(t)=\sum_{k=0}^{+\infty}\, b_k t^k, 
\end{equation}
where the functions $\kappa_1=\kappa_1(t),\ k_1=k_1(t)$ are analytical on $\R$ and their coefficients are connected by the relations
\begin{equation}
\label{3-5}
a_0b_0 = 1,\ \sum_{k=0}^n\Gamma(k+1-\alpha)\Gamma(\alpha+n-k)a_{n-k}b_k = 0,\ n\ge 1.
\end{equation}
The most prominent pair of the Sonine kernels from this class are given by the formulas 
\begin{equation}
\label{Bess}
\kappa(t) = (\sqrt{t})^{\alpha-1}J_{\alpha-1}(2\sqrt{t}),\ 
k(t) = (\sqrt{t})^{-\alpha}I_{-\alpha}(2\sqrt{t}),\ 0<\alpha <1,
\end{equation}
where 
\begin{equation}
\label{Bessf}
J_\nu (t) = \sum_{k=0}^{+\infty} \frac{(-1)^k(t/2)^{2k+\nu}}{k!\Gamma(k+\nu+1)},\ 
I_\nu (t) = \sum_{k=0}^{+\infty} \frac{(t/2)^{2k+\nu}}{k!\Gamma(k+\nu+1)},\ \Re(\nu)>-1,\ t\in \C
\end{equation}
are the Bessel and  the  modified Bessel functions, respectively. 

Later on, the evolution equations with the integro-differential operators of the convolution type (compare to the Sonine solution formula \eqref{Son2})
\begin{equation}
\label{FDR-L} 
(\D_{(k)}\, f)(t) = \frac{d}{dt}\, \int_0^t k(t-\tau)f(\tau)\, d\tau,\ t>0
\end{equation}
were actively studied in the framework of the abstract Volterra integral equations on the Banach spaces (see \cite{Pr} and references therein). For example, in \cite{Cle84}, the case of the operators with the completely  positive kernels $k\in L^1(0,+\infty)$ was treated. The kernels from this class satisfy the condition (compare to the Sonine condition \eqref{Son})
\begin{equation}
\label{Cle_1} 
a\, k(t) \, +\, \int_0^t k(t-\tau)l(\tau)\, d\tau\, = \{1 \},\  t>0,
\end{equation}
where $a \ge 0$ and $l\in L^1(0,+\infty)$ is a non-negative and non-increasing function.
 
However, until recently, no interpretation of these general results in the framework of Fractional Calculus (FC) was suggested. The situation changed with publication of the paper \cite{Koch11} (see also \cite{Han20,Luc21a,LucYam20}). In \cite{Koch11}, Kochubei introduced a class $\mathcal{K}$ of the kernels that satisfy the following conditions:

\vspace{0.1cm}

\noindent
K1) The Laplace transform $\tilde k$ of $k$,
\begin{equation}
\label{Laplace} 
\tilde k(p) = ({\mathcal L}\, k)(p)\ =\ \int_0^{+\infty} k(t)\, e^{-pt}\, dt
\end{equation}
exists for all $p>0$,

\vspace{0.1cm}

\noindent
K2) $\tilde k(p)$ is a Stieltjes function (see \cite{[SSV]} for details regarding the Stieltjes functions),

\vspace{0.1cm}

\noindent
K3) $\tilde k(p) \to 0$ and $p \tilde k(p) \to +\infty$ as $p \to +\infty$,

\vspace{0.1cm}

\noindent
K4)  $\tilde k(p) \to +\infty$ and $p \tilde k(p) \to 0$ as $p \to 0$.

\vspace{0.1cm} 

Using the technique of the complete Bernstein functions, Kochubei investigated   the integro-differential operators in form \eqref{FDR-L} and their Caputo type modifications 
\begin{equation}
\label{FDC}
( _*\D_{(k)}\, f) (t) =  (\D_{(k)}\, f) (t) - f(0)k(t)
\end{equation}
with the kernels from $\mathcal{K}$. In \cite{Koch11}, he showed the inclusion $\mathcal{K}\subset \mathcal{S}$, introduced the corresponding integral operator
\begin{equation}
\label{GFI}
(\I_{(\kappa)}\, f) (t) =  (\kappa\, *\, f)(t) = \int_0^t \kappa(t-\tau)f(\tau)\, d\tau
\end{equation}
and proved validity of the 1st fundamental theorem of FC, i.e., that  the operators \eqref{FDR-L} and \eqref{FDC} are left-inverse to the integral operator \eqref{GFI} on the suitable spaces of functions. 

Moreover, Kochubei treated some basic ordinary and partial fractional differential equations with the time-derivative in form \eqref{FDC} and proved that  the solution to the Cauchy problem for the relaxation equation with the operator \eqref{FDC} and a positive initial condition is completely monotonic and that the fundamental solution to the Cauchy problem for the fractional diffusion equation with the time-derivative in form \eqref{FDC} can be interpreted as a probability density function. These results justified calling the operators \eqref{FDR-L} and \eqref{FDC} the general fractional derivatives (GFDs) in the Riemann-Liouville and Caputo sense, respectively. The integral operator \eqref{GFI} was called the general fractional integral (GFI).

The GFDs \eqref{FDR-L} and \eqref{FDC} with the kernels $k\in \mathcal{K}\subset \mathcal{S}$  possess a series of important properties. However, the conditions K1)-K4)  are very strong (especially the condition K2) and thus in the subsequent publications the operators \eqref{FDR-L} and \eqref{FDC} with the Sonine kernels from some larger classes were considered from the viewpoint of FC and its applications. In  \cite{LucYam16}, a class of the kernels was introduced that ensures validity of a maximum principle for the general time-fractional diffusion equations with the operators of type \eqref{FDC}. Another important class of the Sonine kernels was described in \cite{Han20} in terms of the completely monotone functions. As shown in \cite{Han20}, any singular (unbounded in a neighborhood of the point zero) locally integrable
completely monotone function $\kappa$ is a Sonine kernel and its associate kernel $k$ is also a locally integrable completely monotone function. 

In the recent publications \cite{Luc21a,Luc21b}, the operators \eqref{FDR-L} and \eqref{FDC} with the Sonine kernels from the class $\mathcal{S}_{-1}\subset \mathcal{S}$ that satisfy just some minimal restrictions were studied from the viewpoint of FC. The Sonine kernels $\kappa,\ k\in \mathcal{S}_{-1}$  are continuous on $\R_+$ and possess the integrable singularities of the power function type at the point zero. In particular, in \cite{Luc21a}, the 1st and the 2nd fundamental theorems of FC for the  operators \eqref{FDR-L} and \eqref{FDC} with the kernels $k\in \mathcal{S}_{-1}$ were formulated and proved. In \cite{Luc21b},  an operational calculus of the Mikusi\'nski type for the operators \eqref{FDC} with the Sonine kernels $k\in \mathcal{S}_{-1}$ was constructed and applied for analytical treatment of some initial value problems for the fractional differential equations with these operators. 

It is clear that weakening the Kochubei conditions K1)-K4) on the Sonine kernels from $\mathcal{K}$ leads to giving up of some properties that were derived in \cite{Koch11} for the GFDs \eqref{FDR-L} and \eqref{FDC}. However, it was shown in \cite{Luc21a,Luc21b}  that the  operators \eqref{FDR-L} and \eqref{FDC} with the 
 Sonine kernels $k\in \mathcal{S}_{-1}$  and the corresponding integral operator \eqref{GFI} still satisfy the main properties that the fractional derivatives and integrals should fulfill (see \cite{HL19} and the references therein). Thus, also these operators can be interpreted as the GFDs and GFIs.
 
Another important point concerns the "generalized order" of  the GFDs \eqref{FDR-L} and \eqref{FDC} with the Sonine kernels from the classes mentioned above. While projecting these operators to the conventional Riemann-Liouville and Caputo fractional derivatives (the case of the kernel $k(t) = h_{1-\alpha}(t)$), the derivatives orders are restricted just to the case of $\alpha\in (0,\, 1)$. The reason is that the Sonine condition \eqref{A3} for the power functions $h_\alpha$ and $h_{1-\alpha}$ holds true only in the case $0<\alpha <1$. Moreover, even in the definition of the Caputo type general fractional derivative \eqref{FDC}, only one initial condition is contained that again indicates that the "generalized order" of this operator does not exceed one. 

Because the Riemann-Liouville fractional integral and the Riemann-Liouville and Caputo fractional derivatives are defined for arbitrary  orders $\alpha \ge 0$,  an extension of the GFDs \eqref{FDR-L} and \eqref{FDC} to the case of arbitrary orders is worth for investigation. 

In the recent paper \cite{Luc21a}, the $n$-fold GFIs and GFDs were introduced as an attempt to extend their orders behind the interval $(0,\, 1)$. For example, the two-fold general fractional derivative constructed for the operator \eqref{FDR-L} with the kernel $\kappa(t) = h_{1-\alpha}(t),\ 0<\alpha <1$  is the Riemann-Liouville fractional derivative of the order $2\alpha$:
\begin{equation}
\label{RLn}
(D_{0+}^{2\alpha}\, f)(t) = 
\begin{cases} 
\frac{d^2}{dt^2}(I^{2-2\alpha}_{0+}\, f)(t),& \frac{1}{2} < \alpha <1,\ t>0, \\
 \frac{d}{dt}(I^{1-2\alpha}_{0+}\, f)(t),& 0 < \alpha \le \frac{1}{2},\ t>0. 
 \end{cases}
\end{equation}
Thus, we cannot ensure that the order of this two-fold GFD is always greater than one. Depending on the values of $\alpha$ and $n$, the "generalized order" of the $n$-fold GFD  can be any number from the interval $(0,\, n)$. 

The main objective of this paper is introducing the GFIs and GFDs of arbitrary order in analogy to the Riemann-Liouville fractional integral and the Riemann-Liouville and Caputo fractional derivatives. This is done by a suitable generalization of the Sonine condition \eqref{Son} and by the corresponding  adjustment of the formulas \eqref{FDR-L} and \eqref{FDC} defining the GFDs in the Riemann-Liouville and Caputo senses. 

The rest of the paper is organized as follows. In the 2nd Section, following \cite{Luc21a,Luc21b}, we provide some basic definitions and properties of the GFDs \eqref{FDR-L} and \eqref{FDC} with the Sonine kernels $k\in \mathcal{S}_{-1}$. The 3rd Section contains our main results. First, a suitable generalization of the Sonine condition \eqref{Son} is introduced and some examples of the kernels that satisfy this condition are discussed. Then the GFDs of arbitrary order with these kernels are defined and their properties are studied. The conventional Riemann-Liouville and Caputo fractional derivatives of arbitrary orders are particular cases of these GFDs. Another important example are the integro-differential operators of convolution type with the Bessel and the modified Bessel functions in the kernels. The constructions introduced in this section allow formulation of the fractional differential equations with the GFDs of the generalized order greater than one with several initial conditions.

\section{General Fractional Integrals and Derivatives with the Sonine Kernels}
\label{sec:2}

In this section, we provide some basic definitions and results regarding the GFIs and GFDs with the Sonine 
kernels from the class $\mathcal{S}_{-1}$ introduced in \cite{Luc21a}. For more details, other relevant results, and the proofs see \cite{Luc21a,Luc21b}. 

In what follows, we employ the space of functions $C_{-1}(0,+\infty)$ and its sub-spaces. A family of the spaces $C_\alpha(0,+\infty),\ \alpha \ge -1$ was first introduced in  \cite{Dim66} as follows:
\begin{equation}
\label{2-7}
C_{\alpha}(0,+\infty)\, := \, \{f:\ f(t)=t^{p}f_1(t),\ t>0,\ p>\alpha,\ f_1\in C[0,+\infty)\}.
\end{equation}
Evidently, the spaces $C_{\alpha}(0,+\infty)$ are ordered by inclusion: $\alpha_1 \ge \alpha_2$ $ \Rightarrow$ $C_{\alpha_1}(0,+\infty)\subseteq C_{\alpha_2}(0,+\infty)$ and thus the inclusion $C_{\alpha}(0,+\infty)\subseteq C_{-1}(0,+\infty),\ \alpha \ge -1$ holds true.

In the further discussions, we also use the sub-spaces $C_{-1}^m(0,+\infty),\ m\in \N_0=\N \cup\{0\}$ of the space $C_{-1}(0,+\infty)$ that are defined as follows:
\begin{equation}
\label{2-7-1}
C_{-1}^m(0,+\infty)\, := \, \{f:\ f^{(m)} \in C_{-1}(0,+\infty)\}.
\end{equation}
The spaces $C_{-1}^m(0,+\infty)$ were first introduced and studied in \cite{LucGor99}. In particular, we have the following properties:
\vskip 0.1cm 

\noindent
(1) $C_{-1}^0(0,+\infty)\, \equiv \, C_{-1}(0,+\infty)$.
\vskip 0.1cm 

\noindent
(2) $C_{-1}^m(0,+\infty), \ m\in \N_0$ is a vector space over the field $\R$ (or $\Com$).
\vskip 0.1cm 

\noindent
(3) If $f\in C_{-1}^m(0,+\infty)$ with $m\ge 1$, 
then $f^{(k)}(0+) := \lim\limits_{t\to 0+} f^{(k)}(t) <+\infty,\ 0\le k\le 
m-1$, and the function
$$
\tilde f (t) = 
\begin{cases}f(t), & t>0, \\
f(0+), & t=0
\end{cases}
$$
belongs to the space $C^{m-1}[0,+\infty)$.

\vskip 0.1cm 

\noindent
(4) If $f\in C_{-1}^m(0,+\infty)$ with $m\ge 1$, then $f \in 
C^{m}(0,+\infty)\cap C^{m-1}[0,+\infty)$.
\vskip 0.1cm 

\noindent
(5) For $ m\ge 1$, the following representation holds true:
$$
f\in C_{-1}^m(0,+\infty) \Leftrightarrow
f(t) = (I^m_{0+} \phi)(t) + \sum \limits_{k=0}^{m-1} f^{(k)}(0)\frac{t^k}{k!},\ t\ge 
0,\ \phi \in C_{-1}(0,+\infty).
$$
\vskip 0.1cm 

\noindent
(6) Let $f\in C_{-1}^m(0,+\infty), \ m\in \N_0$,
$f(0)=\dots=f^{(m-1)}(0)=0$ and $g\in C_{-1}^1$. Then,  the 
Laplace convolution $h(t) = (f*g)(t)$
belongs to the space $C_{-1}^{m+1}(0,+\infty)$ and $h(0)=\dots=h^{(m)}(0)=0$.
\vskip 0.1cm 

\noindent
For our aims, we also need another two-parameters family of sub-spaces of $C_{\alpha}(0,+\infty)$ that allows to better control  behavior of the functions at the origin: 
\begin{equation}
\label{subspace}
 C_{\alpha,\beta}(0,+\infty) \, = \, \{f:\ f(t) = t^{p}f_1(t),\ t>0,\ \alpha < p < \beta,\ f_1\in C[0,+\infty)\}.
\end{equation}
In particular, the sub-space $C_{-1,0}(0,+\infty)$ contains the functions continuous on $\R_+$ that possess the integrable singularities of the power function type at the origin. 

% In particular, we mention the well-known formula
% \begin{equation}
% \label{2-13-1}
% \{1\}^n(t) = h_{1}^n(t)= h_{n}(t)=\frac{t^{n-1}}{\Gamma(n)} = \frac{t^{n-1}}{(n-1)!},\ n\in \N.
% \end{equation}

As mentioned in \cite{Sam} (see also \cite{Han20}), any Sonine kernel  has an integrable singularity at the point zero. On the other hand, the kernels of the fractional integrals and derivatives should be singular (\cite {DGGS}). Thus, the fractional integrals and derivatives with the Sonine kernels are worth   being  investigated. In what follows, we consider the GFI \eqref{GFI} and the GFDs \eqref{FDR-L} and \eqref{FDC} of the Riemann-Liouville and Caputo types, respectively, with the Sonine kernels $\kappa$ and $k$ that belong to the sub-space $C_{-1,0}(0,+\infty)$ of the space  $C_{-1}(0,+\infty)$.  

\begin{definition}
\label{dd2}
Let $\kappa,\, k \in C_{-1,0}(0,+\infty)$ be a pair of the Sonine kernels, i.e., the Sonine condition \eqref{Son} be fulfilled. The set of such Sonine kernels    is denoted by $\mathcal{S}_{-1}$:
\begin{equation}
\label{Son_1}
(\kappa,\, k \in \mathcal{S}_{-1} ) \ \Leftrightarrow \ (\kappa,\, k \in C_{-1,0}(0,+\infty))\wedge ((\kappa\, *\, k)(t) \, = \, \{1\}).
\end{equation}
\end{definition}

%\vspace*{-3pt} %%%%%%%%%%%%%%%%

% \section{Fractional Integrals and Derivatives with the Sonine Kernels}
% \label{sec:4}

% \setcounter{section}{4}
% \setcounter{equation}{0}\setcounter{theorem}{0}

%\vspace*{-3pt} %%%%%%%%%%%%%%%%

% In this section, we address the GFIs {GFI} and the GFDs with the Sonine kernels from $\mathcal{S}_{-1}$ on the spaces of functions $C_{-1}^m(0,+\infty), \ m\in \N_0$  and on their sub-spaces. 
% \begin{definition}
% \label{dd3}
% Let $\kappa \in \mathcal{S}_{-1}$. The GFI with the kernel $\kappa$ is defined by the formula
% \begin{equation}
% \label{GFI}
% (\I_{(\kappa)}\, f)(t) := (\kappa\, *\, f)(t) = \int_0^t \kappa(t-\tau)f(\tau)\, d\tau,\ t>0.
% \end{equation}
% \end{definition}

% The GFI \eqref{GFI} with the kernel $\kappa(t) = h_\alpha(t) =\frac{t^{\alpha-1}}{\Gamma(\alpha)}\in \mathcal{S}_{-1},\ 0<\alpha<1$ is the Riemann--Liouville fractional integral \eqref{A4}:
% \begin{equation}
% \label{GFI_1}
% (\I_{(\kappa)}\, f)(t) = (I^\alpha_{0+}\, f)(t),\ t>0.
% \end{equation} 

Several important features of the GFI \eqref{GFI} on the space $C_{-1}(0,+\infty)$ follow from the well-known properties of the Laplace convolution. In particular, we mention the mapping property
\begin{equation}
\label{GFI-map}
\I_{(\kappa)}:\, C_{-1}(0,+\infty)\, \rightarrow C_{-1}(0,+\infty),
\end{equation}
the commutativity law 
\begin{equation}
\label{GFI-com}
\I_{(\kappa_1)}\, \I_{(\kappa_2)} = \I_{(\kappa_2)}\, \I_{(\kappa_1)},\ \kappa_1,\, \kappa_2 \in \mathcal{S}_{-1},
\end{equation}
and the index law
\begin{equation}
\label{GFI-index}
\I_{(\kappa_1)}\, \I_{(\kappa_2)} = \I_{(\kappa_1*\kappa_2)},\ \kappa_1,\, \kappa_2 \in \mathcal{S}_{-1}
\end{equation}
that are valid on the space $C_{-1}(0,+\infty)$.

Let $\kappa \in \mathcal{S}_{-1}$ and $k$ be its associate Sonine kernel. The GFDs of the Riemann-Liouville and the Caputo types associate to the GFI \eqref{GFI} are given by the formulas \eqref{FDR-L} and \eqref{FDC}, respectively.  It is easy to see that the GFD \eqref{FDC} in the Caputo sense can be rewritten as a regularized GFD \eqref{FDR-L} in the Riemann--Liouville sense:
\begin{equation}
\label{GFDC_new}
( _*\D_{(k)}\, f) (t) = (\D_{(k)}\, [f(\cdot)-f(0)]) (t),\ t>0.
\end{equation}
For the functions from $C_{-1}^1(0,\infty)$, the Riemann-Liouville GFD \eqref{FDR-L} can be represented as
\begin{equation}
\label{GFDL-1}
(\D_{(k)}\, f) (t) = (k\, * \, f^\prime)(t) + f(0)k(t),\ t>0,
\end{equation}
that immediately leads to the useful representation 
\vspace{6pt}
\begin{equation}
\label{GFDC_1}
( _*\D_{(k)}\, f) (t) = (k\, * \, f^\prime)(t),\ t>0 
\end{equation}
of the Caputo type GFD \eqref{FDC} that is valid on the space $C_{-1}^1(0,\infty)$. 

In the rest of this section, we formulate the 1st and the 2nd fundamental theorems of FC for the GFDs in the Riemann-Liouville and Caputo senses. 

\begin{theorem}[First Fundamental Theorem for the GFD]
\label{t3}
Let $\kappa \in \mathcal{S}_{-1}$ and $k$ be its associate Sonine kernel. 

Then,  the GFD \eqref{FDR-L} is a left inverse operator to the GFI \eqref{GFI} on the space $C_{-1}(0,+\infty)$: 
\begin{equation}
\label{FTL}
(\D_{(k)}\, \I_{(\kappa)}\, f) (t) = f(t),\ f\in C_{-1}(0,+\infty),\ t>0,
\end{equation}
and the GFD \eqref{FDC} is a left inverse operator to the GFI \eqref{GFI} on the space $C_{-1,k}(0,+\infty)$: 
\begin{equation}
\label{FTC}
( _*\D_{(k)}\, \I_{(\kappa)}\, f) (t) = f(t),\ f\in C_{-1,k}(0,+\infty),\ t>0,
\end{equation}
where $C_{-1,k}(0,+\infty) := \{f:\ f(t)=(\I_{(k)}\, \phi)(t),\ \phi\in C_{-1}(0,+\infty)\}$.
\end{theorem}

As shown in \cite{Luc21a}, the space $C_{-1,k}(0,+\infty)$ can be also characterized as follows:
$$
C_{-1,k}(0,+\infty) = 
 \{f:\ \I_{(\kappa)} f \in C_{-1}^1(0,+\infty)\, \wedge \, (\I_{(\kappa)}\, f)(0) = 0\}.
$$

Now,  we proceed with the 2nd fundamental theorem of FC for the GFDs in the Riemann--Liouville and Caputo senses.

\begin{theorem}[Second  Fundamental Theorem for the GFD]
\label{t4}
Let $\kappa \in \mathcal{S}_{-1}$ and $k$ be its associate Sonine kernel.

Then,  the relations
\begin{equation}
\label{2FTC}
(\I_{(\kappa)}\, _*\D_{(k)}\, f) (t) = f(t)-f(0),\ t>0,
\end{equation}
\begin{equation}
\label{2FTL}
(\I_{(\kappa)}\, \D_{(k)}\, f) (t) = f(t),\ t>0
\end{equation}
hold valid for the functions $f\in C_{-1}^1(0,+\infty)$.
\end{theorem}

In \cite{Luc21a,Luc21b}, the $n$-fold GFIs and GFDs with the Sonine kernels from $\mathcal{S}_{-1}$ were introduced and studies. For details we refer the interested readers to these publications. 

\section{General Fractional Integrals and derivatives of Arbitrary Order}
\label{sec:4}

As already mentioned in Introduction, the "generalized orders" of the GFIs and GFDs introduced so far are restricted to the interval $(0,\, 1)$. The orders of the $n$-fold GFIs and GFDs recently introduced in \cite{Luc21a} belong to the interval $(0,\, n)$. However, it is hardly possible to fix their orders between two neighboring natural numbers as in the case of the conventional Riemann-Liouville and Caputo fractional derivatives and thus to study, say, the fractional oscillator equations or the time-fractional diffusion-wave equations with the GFDs of the order from the interval $(1,\, 2)$.  

In this section, we define the  GFIs and GFDs of arbitrary  orders and study their basic properties. As in the case of the conventional Riemann-Liouville and Caputo fractional derivatives, also for the GFDs we have to distinguish between two completely different cases, namely, between the case of the integer orders and the case of non-integer orders.   In the first case, the conventional Riemann-Liouville and Caputo fractional derivatives are defined as the integer order derivatives while in the second case they are non-local integro-differential operators. Because the conventional Riemann-Liouville and Caputo fractional derivatives are important particular  cases of the GFDs, we have no other choice as to follow the same strategy, namely, to separately define the GFDs of integer orders as the integer orders derivatives and the GFDs of non-integer orders as some integro-differential operators. In  what follows, we focus on the case of the GFDs of non-integer orders (the integer orders GFDs are just the integer order derivatives).   

To introduce the GFIs and the GFDs of arbitrary non-integer orders, we first formulate a condition on their kernels  that generalizes the Sonine condition \eqref{Son}:
\begin{equation}
\label{Luc}
(\kappa \, * \, k)(t) = \{ 1\}^n(t),\ n\in \N,\ t>0,
\end{equation}
where
$$
\{ 1\}^n(t):= (\underbrace{\{ 1\}*\ldots\ * \{ 1\}}_{\mbox{$n$ times}})(t) = h_{n}(t) = \frac{t^{n-1}}{(n-1)!}.
$$
Evidently, the Sonine condition corresponds to the case $n=1$ of the more general condition \eqref{Luc}.

Another important ingredient of our definitions is a set of the kernels that satisfy the condition \eqref{Luc} and belong to the suitable spaces of functions. 

\begin{definition}
\label{d_c}
Let the functions $\kappa$ and $k$ satisfy the condition \eqref{Luc} and the inclusions $\kappa \in C_{-1}(0,+\infty)$, $k\in C_{-1,0}(0,+\infty)$ hold true.

The set of pairs $(\kappa,\, k)$ of such kernels will be denoted by $\mathcal{L}_n$. 
\end{definition}

\begin{remark}
\label{r1}
The set $\mathcal{L}_1$ coincides with the set of the Sonine kernels $\mathcal{S}_{-1}$ that we discussed in the previous section (see Definition \ref{dd2}). Indeed, in this case, the kernel $\kappa \in C_{-1}(0,+\infty)$ is a Sonine kernel and therefore it has an integrable  singularity at the point zero. Thus, it belongs to the subspace $C_{-1,0}(0,+\infty)$ as required in Definition \ref{dd2}. 
\end{remark}

\begin{remark}
\label{r2}
For $n>1$, Definition  \ref{d_c} is not symmetrical with respect to the kernels $\kappa$ and $k$  because of the non-symmetrical inclusions  $\kappa \in C_{-1}(0,+\infty)$ and $k\in C_{-1,0}(0,+\infty)$ (in the case $n=1$, Definition \ref{dd2} is symmetrical and one can interchange the kernels $\kappa$ and $k$). 

However, the same statement is valid for the kernel $\kappa(t) = h_{\alpha}(t),\ \alpha >0$ of the Riemann-Liouville integral $I_{0+}^\alpha$  and the kernel $k(t) = h_{n-\alpha}(t)$ of the Riemann-Liouville and Caputo fractional derivatives of order $\alpha,\ n-1<\alpha<n,\ n \in \N$ defined as follows:
\begin{equation}
\label{RLD}
(D^\alpha_{0+}\, f)(t) := \frac{d^n}{dt^n}(I^{n-\alpha}_{0+} \, f)(t),\ t>0,
\end{equation}
\begin{equation}
\label{CD}
( _*D^\alpha_{0+}\, f)(t) := \left(D^\alpha_{0+}\,\left(  f(\cdot) - \sum_{j=0}^{n-1}f^{(j)}(0)h_{j+1}(\cdot)\right)\right)(t),\ t>0, 
\end{equation}
\vspace{6pt}
$I^\alpha_{0+}$ being the Riemann-Liouville fractional integral of order $\alpha$:
\begin{equation}
\label{RLI}
 (I_{0+}^\alpha\, f)(t) := \frac{1}{\Gamma(\alpha)}\, \int_0^t (t-\tau)^{\alpha-1}\, f(\tau)\, d\tau, \ t>0,\ \alpha >0.
\end{equation}
The trick in defining the integer orders  Riemann-Liouville and Caputo fractional derivatives consists in a separate definition of the Riemann-Liouville fractional integral of the order $\alpha = 0$:
\begin{equation}
\label{RLI-0}
 (I_{0+}^0\, f)(t) := f(t).
\end{equation}
Of course, the definition \eqref{RLI-0} is not arbitrary and justified i.a. by the formula
\begin{equation}
\label{RLzerolim}
\lim_{\alpha \to 0+}  \|I_{0+}^\alpha\, f)(t) - f(t)\|_{L_1(0,T)}= 0
\end{equation}
that is valid for $f\in L^1(0,T)$ in every Lebesgue point of $f$, i.e., almost everywhere on the interval $(0,\, T),\ T>0$ (see e.g. \cite{Samko}). 
\end{remark}

\begin{example}
\label{ex1}
The kernels $\kappa(t) = h_{\alpha}(t),\ \alpha >0$ and $k(t) = h_{n-\alpha}(t), \ n-1<\alpha<n,\ n \in \N$ provide a first example of the kernels from $\mathcal{L}_n$. Please note that the power functions $h_{\alpha}$ and $h_{n-\alpha}$ build a pair of the Sonine kernels only in the case $n=1$, i.e., only in the case when the fractional derivatives orders are less than one. 
\end{example}

Because both the Sonine condition \eqref{Son} and its generalization \eqref{Luc} contain the Laplace convolution of two kernels, it is very natural to transform them into the Laplace domain.  Providing the Laplace transforms $\tilde \kappa,\ \tilde k$ of the functions $\kappa$ and $k$   exist, the convolution theorem for the Laplace transform leads to the relation
\begin{equation}
\label{Son-Lap} 
\tilde \kappa(p) \cdot \tilde k(p) = \frac{1}{p},\ \Re(p)>p_{\kappa,k} \in \R
\end{equation}
for the Laplace transforms of the Sonine kernels and to a more general relation
\begin{equation}
\label{Luc-Lap} 
\tilde \kappa(p) \cdot \tilde k(p) = \frac{1}{p^{n}},\ \Re(p)>p_{\kappa,k} \in \R,\ n\in \N
\end{equation} 
for the kernels from the set $\mathcal{L}_n$ introduced in Definition \ref{d_c}. 

\begin{example}
\label{ex2}
The formula \eqref{Luc-Lap} along with the reference books \cite{PBM1}, \cite{PBM2} for the direct and inverse Laplace transforms, respectively, can be used to deduce other nontrivial examples of the kernels from $\mathcal{L}_n$.  For instance, we employ the Laplace transform formulas (see \cite{PBM1})
$$ 
\left({\mathcal L}\ t^{\nu/2}J_{\nu}(2\sqrt{t})\right)(p)\, = \, p^{-\nu -1}\, \exp(-1/p),\ \Re(\nu)>-1,\ \Re(p) >0,
$$
$$ 
\left({\mathcal L}\ t^{\nu/2}I_{\nu}(2\sqrt{t})\right)(p)\, = \, p^{-\nu -1}\, \exp(1/p),\ \Re(\nu)>-1,\ \Re(p) >0
$$
for the Bessel function $J_\nu$ and the modified Bessel function $I_\nu$ defined by the power series \eqref{Bessf} to introduce the kernels
\begin{equation}
\label{Luc-Bes} 
\kappa(t) = t^{\nu/2}\, J_{\nu}(2\sqrt{t}),\ \ k(t) = t^{n/2-\nu/2-1}\, I_{n-\nu-2}(2\sqrt{t}),\ n-2 < \nu < n-1,\ n\in \N.
\end{equation} 
These kernels satisfy the condition \eqref{Luc-Lap}. Moreover, for $n-2 < \nu < n-1,\ n\in \N$ the inclusions  $\kappa \in C_{-1}(0,+\infty)$ and $k\in C_{-1,0}(0,+\infty)$ hold true and thus the pair of the kernels $(\kappa,\ k)$ given by \eqref{Luc-Bes} is from $\mathcal{L}_n$. 
\end{example}

Now let us consider a pair of the Sonine kernels $(\kappa,\ k)$ from $\mathcal{L}_1$ (in \cite{Han20,Luc21a,Luc21b,Sam,Son} and other related publications, many pairs of such kernels were presented). There are at least two reasonable possibilities to construct a pair $(\kappa_n\, k_n)$ of the kernels from $\mathcal{L}_n,\ n>1$ based on the Sonine kernels $\kappa,\, k$ from $\mathcal{L}_1$. 

The first strategy consists in  building the kernels $\kappa_n = \kappa^n$ and $k_n=k^n$. 
Evidently, the kernels $\kappa_n$ and $k_n$ satisfy  the relation \eqref{Luc} because $\kappa$ and $k$ are the Sonine kernels:
\begin{equation}
\label{Son_n} 
(\kappa_n\, *\, k_n)(t) = (\kappa^n\, *\, k^n)(t) = (\kappa\, *\, k)^n(t) = \{ 1\}^n(t).
\end{equation} 
However, the pair $(\kappa_n,\ k_n)$ does not always belong to the set $\mathcal{L}_n$. This is the case only under an additional condition, namely, only when the inclusion $k^n\in C_{-1,0}(0,+\infty)$ holds true (of course, $\kappa^n \in C_{-1}(0,+\infty)$ for any $n\in \N$). And this is a very strong and restrictive condition.   For example, in the case of the Riemann-Liouville fractional integral $I_{0+}^\alpha$ with the kernel $\kappa(t)=h_\alpha(t), 0<\alpha<1$ and the Riemann-Liouville fractional derivative $D_{0+}^\alpha$ with the kernel $k(t)=h_{1-\alpha}$, the kernel $k^n$ takes the form $k^n(t) = h_{n(1-\alpha)}(t)$. It belongs to the space $C_{-1,0}(0,+\infty)$ only under the condition $0<n(1-\alpha)<1$, i.e., if $1- \frac{1}{n} <\alpha < 1$ that is very restrictive. Moreover, the example of the kernels \eqref{Luc-Bes} shows that not any pair of the kernels from $\mathcal{L}_n$ can be represented in  form $(\kappa^n,\ k^n)$ with the kernels 
$(\kappa,\ k)\in \mathcal{L}_1$. 

Another and even more general and important possibility  for construction a pair $(\kappa_n,\, k_n)$ of the kernels from $\mathcal{L}_n,\ n>1$ based on the Sonine kernels $\kappa,\, k$ from $\mathcal{L}_1$ is presented in the following theorem:

\begin{theorem}
\label{tkernel}
Let $(\kappa,\ k)$ be a pair of the Sonine kernels from $\mathcal{L}_1$.  

Then the  pair $(\kappa_n,\ k_n)$ of the kernels given by the formula
\begin{equation}
\label{Son_Luc_n}
\kappa_n(t) = (\{1\}^{n-1}\, *\, \kappa)(t),\ \ k_n(t)=k(t)
\end{equation}
belongs to the set $\mathcal{L}_n$.
\end{theorem}

\begin{proof}
First we check that the kernels \eqref{Son_Luc_n} satisfy the condition \eqref{Luc}:
\begin{equation}
\label{Son_n_1} 
(\kappa_n\, *\, k_n)(t) = (\{1\}^{n-1}\, *\, \kappa\, *\, k)(t)  =
(\{1\}^{n-1}\, *\,\{ 1\})(t) = \{ 1\}^n(t).
\end{equation} 
Moreover, because of the inclusions $\kappa,\, k\in \mathcal{L}_1$, the inclusions  $\kappa_n \in C_{-1}(0,+\infty)$ and  $k_n = k \in C_{-1,0}(0,+\infty)$ are satisfied and thus the kernels $\kappa_n$ and $k_n$ defined by \eqref{Son_Luc_n} belong to the set $\mathcal{L}_n$.
\end{proof}

In the rest of this section, we introduce the general fractional integrals and derivatives of arbitrary (non-integer) order and discuss their basic properties and examples.   

\begin{definition}
\label{dao}
Let  $(\kappa,\ k)$ be a pair of the kernels from $\mathcal{L}_n$. The GFI with the kernel $\kappa$ is specified by the standard formula
\begin{equation}
\label{GFIn}
(\I_{(\kappa)}\, f)(t) :=  \int_0^t \kappa(t-\tau)f(\tau)\, 
d\tau,\ t>0,
\end{equation}
whereas the GFDs of the Riemann-Liouville and Caputo types with the kernel $k$ are defined as follows:
\begin{equation}
\label{FDR-Ln} 
(\D_{(k)}\, f)(t) := \frac{d^n}{dt^n}\, \int_0^t k(t-\tau)f(\tau)\, d\tau,\ t>0,
\end{equation}
\begin{equation}
\label{FDCn}
( _*\D_{(k)}\, f)(t) := \left(\D_{(k)}\,\left(  f(\cdot) - \sum_{j=0}^{n-1}f^{(j)}(0)h_{j+1}(\cdot)\right)\right)(t),\ t>0. 
\end{equation}
\end{definition}

\begin{example}
\label{ex3}
Evidently, the  GFI \eqref{GFIn} with the kernel $\kappa(t) = h_\alpha(t),\ \alpha>0$ is reduced to the Riemann-Liouville fractional integral \eqref{RLI} and the Riemann-Liouville and Caputo fractional derivatives of the order $\alpha,\ n-1<\alpha<n,\ n \in \N$ defined by \eqref{RLD} and \eqref{CD}, respectively, are particular cases of the GFDs \eqref{FDR-Ln} and \eqref{FDCn} with the kernel $k(t) = h_{n-\alpha}(t)$. As mentioned in Example \ref{ex1}, the inclusion $(h_\alpha,\ h_{n-\alpha})\in \mathcal{L}_n$ holds valid if and only if $n-1<\alpha<n,\ n \in \N$. 
\end{example}

It is worth mentioning that the Riemann-Liouville fractional integral \eqref{RLI} and the Riemann-Liouville and Caputo fractional derivatives of an arbitrary order $\alpha,\ n-1<\alpha<n,\ n \in \N$ can be introduced based on the Sonine pair $\kappa =h_\beta$, $k=k_{1-\beta}$, $0<\beta < 1$ and using the construction \eqref{Son_Luc_n} presented in Theorem \ref{tkernel}. Indeed, in this case we have the relations 
\begin{equation}
\label{Son_Luc_n_RL}
\kappa_n(t) = (\{1\}^{n-1}\, *\, \kappa)(t) = (\{1\}^{n-1}\, *\, h_{\beta})(t) = h_{n-1+\beta}(t),\ k_n(t)=k(t)=h_{1-\beta}(t).
\end{equation}
Thus, the GFI \eqref{GFIn} and the GFDs \eqref{FDR-Ln} and \eqref{FDCn} with the kernels $(\kappa_n,\ k_n)\in \mathcal{L}_n$ take the form
\begin{equation}
\label{GFIn_RL}
(\I_{(\kappa)}\, f)(t) =  (h_{n-1+\beta}\, *\, f) (t)  = (I_{0+}^{n-1+\beta}\, f)(t),\ t>0,
\end{equation}
\begin{equation}
\label{FDR-Ln_RL} 
(\D_{(k)}\, f)(t) = \frac{d^n}{dt^n}\, (h_{1-\beta}\, *\, f)(t) = \frac{d^n}{dt^n}\, (I_{0+}^{1-\beta}\, f)(t),\ t>0,
\end{equation}
\begin{equation}
\label{FDCn_RL}
( _*\D_{(k)}\, f)(t) = \frac{d^n}{dt^n}\, \left(I_{0+}^{1-\beta}\, \left(  f(\cdot) - \sum_{j=0}^{n-1}f^{(j)}(0)h_{j+1}(\cdot)\right)\right)(t),\ t>0. 
\end{equation}
Now we introduce a new variable $\alpha:=n-1+\beta$. Then $1-\beta = n-\alpha$ and the inequalities $n-1<\alpha <n$ are fulfilled because of the condition $0<\beta < 1$. Thus, the operator \eqref{GFIn_RL} is the Riemann-Liouville fractional integral \eqref{RLI} of the order $\alpha$ and the operators \eqref{FDR-Ln_RL} and \eqref{FDCn_RL} coincide with the Riemann-Liouville and Caputo fractional derivatives of the order $\alpha,\ n-1<\alpha <n, \ n\in \N$. 

\begin{example}
\label{ex4}
Another interesting and nontrivial particular case of the GFI \eqref{GFIn} and the GFDs \eqref{FDR-Ln} and \eqref{FDCn} is constructed for the pair $(\kappa,\ k)\in \mathcal{L}_n$ of the kernels defined by the formula \eqref{Luc-Bes} with  $n-2 < \nu < n-1,\ n\in \N$: 
\begin{equation}
\label{GFIn_B}
(\I_{(\kappa)}\, f)(t) =  \int_0^t (t-\tau)^{\nu/2}\, J_{\nu}(2\sqrt{t-\tau})f(\tau)\, 
d\tau,\ t>0,
\end{equation}
\begin{equation}
\label{FDR-Ln_B} 
(\D_{(k)}\, f)(t) = \frac{d^n}{dt^n}\, \int_0^t (t-\tau)^{n/2-\nu/2-1}I_{n-\nu-2}(2\sqrt{t-\tau})f(\tau)\, d\tau,\ t>0,
\end{equation}
\begin{equation}
\label{FDCn_B}
( _*\D_{(k)}\, f)(t) := \left(\D_{(k)}\,\left(  f(\cdot) - \sum_{j=0}^{n-1}f^{(j)}(0)h_{j+1}(\cdot)\right)\right)(t),\ t>0. 
\end{equation}
\end{example}

It is worth mentioning that the Caputo type GFD \eqref{FDCn} can be represented in a slightly different form:
$$
( _*\D_{(k)}\, f)(t) = \left(\D_{(k)}\,\left(  f(\cdot) - \sum_{j=0}^{n-1}f^{(j)}(0)h_{j+1}(\cdot)\right)\right)(t) = 
$$
$$
(\D_{(k)}\, f)(t) - \sum_{j=0}^{n-1} f^{(j)}(0)(\D_{(k)}\, h_{j+1})(t) =
(\D_{(k)}\, f)(t) - \sum_{j=0}^{n-1} f^{(j)}(0)\, \frac{d^n}{dt^n} (k\, *\, h_{j+1})(t) = 
$$
\begin{equation}
\label{FDCn_1}
(\D_{(k)}\, f)(t) - \sum_{j=0}^{n-1} f^{(j)}(0)\, \frac{d^n}{dt^n} (I_{0+}^{j+1}\, k)(t) =
(\D_{(k)}\, f) (t) - \sum_{j=0}^{n-1}f^{(j)}(0) 
\frac{d^{n-j-1}}{dt^{n-j-1}} k(t),\ t>0.
\end{equation}

As to the basic properties  of the GFI \eqref{GFIn} of arbitrary order on  $C_{-1}(0,+\infty)$, they follow from  the well-known properties of the Laplace convolution (compare to the properties of the GFI \eqref{GFI} of the order less than one):  
\begin{equation}
\label{GFI-map_1}
\I_{(\kappa)}:\, C_{-1}(0,+\infty)\, \rightarrow C_{-1}(0,+\infty)\ \mbox{(mapping property)},
\end{equation}
\begin{equation}
\label{GFI-com_1}
\I_{(\kappa_1)}\, \I_{(\kappa_2)} = \I_{(\kappa_2)}\, \I_{(\kappa_1)} \ 
\mbox{(commutativity law)}, 
\end{equation}
\begin{equation}
\label{GFI-index_1}
\I_{(\kappa_1)}\, \I_{(\kappa_2)} = \I_{(\kappa_1*\kappa_2)}\ \mbox{(index law)}.
\end{equation}

To justify the denotation GFIs and GFDs, in the rest of this section, we formulate and prove the two fundamental theorems of FC for the GFDs \eqref{FDR-Ln} and \eqref{FDCn} of the Riemann-Liouville and Caputo types.  

\begin{theorem}[First Fundamental Theorem for the GFD of arbitrary order]
\label{t3_n}

Let  $(\kappa,\ k)$ be a pair of the kernels from $\mathcal{L}_n$. 

Then,  the GFD \eqref{FDR-Ln} is a left inverse operator to the GFI \eqref{GFIn} on the space $C_{-1}(0,+\infty)$: 
\begin{equation}
\label{FTLn}
(\D_{(k)}\, \I_{(\kappa)}\, f) (t) = f(t),\ f\in C_{-1}(0,+\infty),\ t>0,
\end{equation}
and the GFD \eqref{FDCn} is a left inverse operator to the GFI \eqref{GFIn} on the space $C_{-1,k}(0,+\infty)$: 
\begin{equation}
\label{FTCn}
( _*\D_{(k)}\, \I_{(\kappa)}\, f) (t) = f(t),\ f\in C_{-1,k}(0,+\infty),\ t>0,
\end{equation}
where the space $C_{-1,k}(0,+\infty)$ is defined as in Theorem \ref{t3}.
\end{theorem}
\begin{proof}
We start with a proof of the formula \eqref{FTLn}:
$$
(\D_{(k)}\, \I_{(\kappa)}\, f) (t) = \frac{d^n}{dt^n}\,(k\, *\, (\kappa\, *\, f))(t) = 
\frac{d^n}{dt^n}\,((k\, *\, \kappa)\, *\, f)(t) = 
$$
$$
\frac{d^n}{dt^n}\,(\{1\}^n\, *\, f)(t) = 
\frac{d^n}{dt^n}\,(I_{0+}^n\, f)(t) = f(t).
$$
A function $f\in C_{-1,k}(0,+\infty)$ can be represented in the form $f(t) = (\I_{(k)}\, \phi)(t),\ \phi\in C_{-1}(0,+\infty)$ and thus the following chain of equations is valid:
$$
(\I_{(\kappa)}\, f)(t) =  (\I_{(\kappa)}\, \I_{(k)}\, \phi)(t) = ((\kappa\, *\, k)\, *\, f)(t) = (\{1\}^n\, \phi)(t) = (I_{0+}^n\, \phi)(t).
$$
The last relation implicates the inclusion $\I_{(\kappa)}\, f \in C_{-1}^n(0,+\infty)$ and the relations
\begin{equation}
\label{hilf}
\frac{d^j}{dt^j}\, (\I_{(\kappa)}\, f)(t)\Bigl|_{t=0} = (I_{0+}^{n-j}\, \phi)(t)\Bigl|_{t=0} = 0,\ j=0,\dots,n-1. \Bigr.\Bigr.
\end{equation}
To derive the formula \eqref{FTCn}, we employ the representation \eqref{FDCn_1} of the GFD of the Caputo type, the formula \eqref{hilf}, and the relation \eqref{FTLn} that we already proved:
$$
( _*\D_{(k)}\, \I_{(\kappa)}\, f) (t) = (\D_{(k)}\, \I_{(\kappa)}\, f) (t) - \sum_{j=0}^{n-1} \frac{d^j}{dt^j}\, (\I_{(\kappa)}\, f)(t)\Bigl|_{t=0} \frac{d^{n-j-1}}{dt^{n-j-1}} k(t)
= f(t). \Bigr.
$$
\end{proof}

\begin{theorem}[Second  Fundamental Theorem for the GFD of arbitrary order]
\label{t4_n}
Let  $(\kappa,\ k)$ be a pair of the kernels from $\mathcal{L}_n$. 

Then,  the relation
\begin{equation}
\label{2FTCn}
(\I_{(\kappa)}\, _*\D_{(k)}\, f) (t) = f(t) - \sum_{j=0}^{n-1} f^{(j)}(0)\, h_{j+1}(t)
\end{equation}
holds true on the space $C_{-1}^n(0,+\infty)$ and the formula
\begin{equation}
\label{2FTLn}
(\I_{(\kappa)}\, \D_{(k)}\, f) (t) = f(t),\ t>0
\end{equation}
is valid for the functions $f\in C_{-1,\kappa}^n(0,+\infty)$.
\end{theorem}
\begin{proof}
As already mentioned in Section \ref{sec:2}, any function $f$ from $C_{-1}^n(0,+\infty)$ can be represented as follows (see \cite{LucGor99}):
\begin{equation}
\label{frep}
f(t) = (I^n_{0+} \phi)(t) + \sum_{j=0}^{n-1} f^{(j)}(0)\, h_{j+1}(t),\ t\ge 0,\ \phi \in C_{-1}(0,+\infty).
\end{equation}
Then we employ this representation and the formula \eqref{FDCn} and arrive at the following chain of relations:
$$
( _*\D_{(k)}\, f) (t) = ( \D_{(k)}\, \left(f(\cdot) - \sum_{j=0}^{n-1} f^{(j)}(0)\, h_{j+1}(\cdot)\right) (t) = ( \D_{(k)}\, I^n_{0+} \phi)(t) = 
$$
$$
\frac{d^n}{dt^n}\, (k\, *\ \{1\}^n\, *\, \phi)(t) = \frac{d^n}{dt^n}\, (\{1\}^n\, *\,(k\, *\, \phi))(t) = (k\, *\, \phi)(t).
$$
Finally, we take into account the representation \eqref{frep} and get the formula \eqref{2FTCn}:
$$
(\I_{(\kappa)}\, _*\D_{(k)}\, f) (t) = (\I_{(\kappa)}\, (k\, *\, \phi)) (t) = 
((\kappa\, *\, k)\, *\, \phi)(t) =
$$
$$
(\{ 1\}^n \, *\, \phi)(t) = (I^n_{0+} \phi)(t) = f(t) - \sum_{j=0}^{n-1} f^{(j)}(0)\, h_{j+1}(t).
$$
To prove the formula \eqref{2FTLn}, we first mention that a function $f\in C_{-1,\kappa}(0,+\infty)$ can be represented in the form $f(t) = (\I_{(\kappa)}\, \phi)(t),\ \phi\in C_{-1}(0,+\infty)$ and thus the following chain of equations is valid:
$$
(\I_{(\kappa)}\, \D_{(k)}\, f) (t) = (\I_{(\kappa)}\, \frac{d^n}{dt^n}\, (k\, *\, f) (t) = (\I_{(\kappa)}\, \frac{d^n}{dt^n}\, (k\, *\, (\kappa \, *\, \phi)) (t)=
$$
$$
(\I_{(\kappa)}\, \frac{d^n}{dt^n}\, (\{1\}^n *\, \phi)) (t)= (\I_{(\kappa)}\, \phi) (t)= f(t).
$$
\end{proof}

In conclusion, we emphasize once again the result from Theorem 
\ref{tkernel} and its implications on the definitions of the GFIs and the GFDs of arbitrary order. If $(\kappa,\ k)$ is a pair of the Sonine kernels from $\mathcal{L}_1$,  the  pair $(\kappa_n,\ k_n)$ of the kernels given by the formula \eqref{Son_Luc_n} belongs to the set $\mathcal{L}_n,\ n>1$. The GFI \eqref{GFIn} with the kernel $\kappa_n = (\{1\}^{n-1}\, *\, \kappa)(t)$ takes the form 
\begin{equation}
\label{GFIn_n}
(\I_{(\kappa_n)}\, f)(t) = (I_{0+}^{n-1}\, \I_{(\kappa)}\, f)(t),\ t>0,
\end{equation}
whereas the GFDs of the Riemann-Liouville and Caputo types with the kernel $k_n = k$ can be represented as follows:
\begin{equation}
\label{FDR-Ln_n} 
(\D_{(k_n)}\, f)(t) = \frac{d^n}{dt^n}\, (\I_{(k)}\, f)(t),\ t>0,
\end{equation}
\begin{equation}
\label{FDCn_n}
( _*\D_{(k_n)}\, f)(t) = \frac{d^n}{dt^n} \left( \I_{(k)}\,  \left( f(\cdot) - \sum_{j=0}^{n-1}f^{(j)}(0)h_{j+1}(\cdot)\right)\right)(t),\ t>0. 
\end{equation}
As we see, these constructions are completely anological to the definitons of the Riemann-Liouville fractional integral and the Riemann-Liouville and Caputo fractional derivatives of arbitrary order.


\begin{thebibliography}{99}
 \normalsize
%\end{paracol}
%\reftitle{References}
%\begin{thebibliography}{999}
\bibitem{Abel1}
N.H. Abel, Opl\"{o}sning af et par opgaver ved hjelp af bestemte integraler,
{\it Magazin for Naturvidenskaberne}, Aargang I, Bind 2, Christiania, 1823.

\bibitem{Abel2}
N.H. Abel, Aufl\"{o}sung einer mechanischen Aufgabe, {\it Journal f\"{u}r die Reine und Angewandte Mathematik} {\bf 1} (1826), 153--157.


% \bibitem{BasLuc95}
% M.A. Al-Bassam, Yu. Luchko, On generalized
% fractional calculus and its application to the solution of
% integro-differential equations. \emph{Journal of Fractional
% Calculus} \textbf{7} (1995), 69--88.

\bibitem{Cle84}
Ph. Cl\'ement,  On abstract Volterra equations in Banach spaces with completely positive kernels. In: F. Kappel, W. Schappacher (Eds), \emph{  Lecture Notes in Math. Vol. 1076},  Springer, Berlin (1984), pp.~32--40.

\bibitem{DGGS}
K. Diethelm, R. Garrappa, A. Giusti, M. Stynes, Why fractional derivatives with nonsingular kernels should not be used. \emph{ Fract. Calc. Appl. Anal.} \textbf{23} (2020), 610--634.

\bibitem{Dim66}
I.H. Dimovski, Operational  calculus  for  a  class  of
differential operators.  \emph{ Compt.  rend.  Acad.  bulg.  Sci.} \textbf{19} (1966), 1111--1114.

% \bibitem{Han}
% L.A-M. Hanna, Yu. Luchko, 
% Operational calculus for the Caputo-type fractional Erd\'elyi-Kober  derivative and its applications. \emph{Integral Transforms and Special Functions} \textbf{25} (2014), 359--373. 

\bibitem{Han20}
A. Hanyga,  A comment on a controversial issue: A Generalized Fractional Derivative cannot have a regular kernel. \emph{ Fract. Calc. Anal. Appl.} \textbf{23} (2020), 211--223.

\bibitem{HL19}
R. Hilfer, Yu. Luchko,
Desiderata for Fractional Derivatives and Integrals. \emph{ Mathematics } \textbf{7(2)} (2019), 149.

% \bibitem{JK17}
% J. Janno, K.  Kasemets,  Identification of a kernel in an evolutionary integral equation occurring in subdiffusion. \emph{ J. Inverse Ill-Posed Probl.} \emph{25} (2017), 777--798.

% \bibitem{KJ19_1}
% N. Kinash, J. Janno,  Inverse problems for a generalized subdiffusion equation with final overdetermination. \emph{  Math. Model. Anal.} \emph{ 24} (2019), 236--262.

% \bibitem{KJ19}
% N. Kinash, J. Janno, 
% An Inverse Problem for a Generalized Fractional Derivative with an Application in Reconstruction of Time- and Space-Dependent Sources in Fractional Diffusion and Wave Equations. \emph{  Mathematics} \emph{7(12)} (2019), 1138.

\bibitem{Koch11}
A.N. Kochubei,  General fractional calculus, evolution equations, and renewal processes. \emph{Integr. Equa. Operator Theory} \textbf{71} (2011), 583--600.

% \bibitem{Koch19_1}
% A.N. Kochubei, General fractional calculus. In: A. Kochubei, Yu. Luchko (Eds.), \emph{Handbook of Fractional Calculus with Applications. 
% Volume 1: Basic Theory}, De Gruyter, Berlin (2019), pp.~111--126.

% \bibitem{Koch19_2}
% A.N. Kochubei,  Equations with general fractional time derivatives. Cauchy problem. In: A. Kochubei, Yu. Luchko (Eds.), \emph{Handbook of Fractional Calculus with Applications. 
% Vol.2: Fractional Differential Equations}, De Gruyter, Berlin (2019), pp.~223--234.

% \bibitem{KK}
% A.N. Kochubei,  Yu. Kondratiev,  Growth Equation of the General Fractional Calculus. \emph{Mathematics} \textbf{7(7)} (2019), 615. 

% \bibitem{Luc93}
% Yu. Luchko,
%  \emph{Theory of the integral transformations with the Fox H-function as a kernel and some of its applications including operational calculus}. PhD. thesis, Belorussian State University, Minsk (1993). 

% \bibitem{Luc99}
% Yu. Luchko,  
% Operational method in fractional calculus. \emph{Fract. Calc. Appl. Anal.}  \textbf{2}
% (1999), 463--489. 

% \bibitem{Luc20}
% Yu. Luchko, Fractional derivatives and the fundamental theorem of Fractional Calculus. \emph{ Fract. Calc. Appl. Anal.} \textbf{23} (2020), 939--966.

\bibitem{Luc21a}
Yu. Luchko,
General Fractional Integrals and Derivatives with the Sonine Kernels. \emph{Mathematics}  \textbf{2021}, 9(6), 594.

\bibitem{Luc21b}
Yu. Luchko, Operational Calculus for the general fractional derivatives with the Sonine kernels. \emph{arXiv:2103.00549v2, Classical Analysis and ODEs (math.CA)} \textbf{2021}. 

\bibitem{LucGor99}
Yu. Luchko, R. Gorenflo, An operational method for
solving fractional differential equations. \emph{Acta Mathematica
Vietnamica}  \textbf{24} (1999), 207--234. 

% \bibitem{LucSri95}
% Yu. Luchko, H.M. Srivastava,  
% The exact solution of certain differential equations of fractional
% order by using operational calculus. \emph{ Comput. Math. Appl.} \textbf{29}
% (1995), 73--85.

% \bibitem{LucYak94}
% Yu. Luchko, S. Yakubovich, An operational
% method for solving some classes of integro-differential equations.
% \emph{ Differential Equations} \textbf{30} (1994), 247--256.

\bibitem{LucYam16}
Yu. Luchko,  M. Yamamoto, General time-fractional diffusion equation: some uniqueness and existence results for the initial-boundary-value problems. \emph{ Fract. Calc. Appl. Anal.}  \textbf{19} (2016), 675--695.

\bibitem{LucYam20}
Yu. Luchko,  M. Yamamoto,
The General Fractional Derivative and Related Fractional Differential Equations.  \emph{ Mathematics} \textbf{ 8(12)} (2020), 2115.

\bibitem{Pod}
I. Podlubny, R.L. Magin, I. Trymorush, Niels Henrik Abel and the birth of Fractional Calculus,
\emph{  Fract. Calc. Appl. Anal.} {\bf 20} (2017), 1068--1075.

\bibitem{PBM1}
A.P. Prudnikov, Yu.A. Brychkov, O.I. Marichev, 
\emph{ Integrals and Series: Direct Laplace transforms, Vol. 4}, Gordon \& Breach, New York (1992).

\bibitem{PBM2}
A.P. Prudnikov, Yu.A. Brychkov, O.I. Marichev, 
\emph{ Integrals and Series: Inverse Laplace Transforms, Vol. 5},  Gordon \& Breach, New York (1992).

\bibitem{Pr}
J. Pr\"uss, \emph{Evolutionary Integral Equations and Applications}, Springer, Basel (1993).

% \bibitem{Rub}
% B.S. Rubin,  G.F. Volodarskaja, An imbedding theorem for convolutions on
% a finite interval and its application to integral equations of the first kind.
% \emph {Soviet Math. Dokl.} \textbf{20} (1979), 234--239.

\bibitem{Sam}
S.G. Samko, R.P. Cardoso, Integral equations of the first kind of Sonine
type. \emph {Intern. J. Math. Sci.} \textbf{57} (2003), 3609--3632.

\bibitem{Samko}
S.G. Samko, A.A. Kilbas, O.I. Marichev,   
\emph {Fractional Integrals and Derivatives. Theory and Applications}, Gordon and Breach, 
New York (1993).

\bibitem{[SSV]}
R.L. Schilling, R. Song, Z. Vondracek, \emph { Bernstein
Functions. Theory and Application}, De Gruyter, Berlin (2010).

% \bibitem{Sin18}
% Ch.-S. Sin,  Well-posedness of general Caputo-type fractional differential equations. 
% \emph{Fract. Calc. Appl. Anal.} \textbf{21} (2018), 819--832.

\bibitem{Son}
N. Sonine, Sur la g\'en\'eralisation d'une formule d'Abel. \emph{ Acta Math.} \textbf{4} (1884), 171--176.

% \bibitem{Wick}
% J. Wick, \"Uber eine Integralgleichung vom Abelschen Typ. \emph{ Z. Angew. Math. Mech.} 
% \textbf{48} (1968), T39--T41.

% \bibitem{YakLuc94}
% S. Yakubovich, Yu. Luchko, 
% \emph{The hypergeometric approach to integral transforms and convolutions},  Kluwer Acad. Publ., Dordrecht (1994).

% \bibitem{Zac08}
% R. Zacher, 
% Boundedness of weak solutions to evolutionary partial integro-differential
% equations with discontinuous coefficients. \emph{ J. Math. Anal. Appl.} \textbf{ 348} (2008), 137--149.

\end{thebibliography}
\end{document}